\documentclass[11pt]{article}

\usepackage{amssymb, amsmath, amsthm, graphicx}
\usepackage[left=1in,top=1in,right=1in]{geometry}

\usepackage{subfigure}

      \theoremstyle{plain}
      \newtheorem{theorem}{Theorem}[section]
      \newtheorem{lemma}[theorem]{Lemma}
            
            \newtheorem{observation}[theorem]{Observation}

      \theoremstyle{definition}

      \theoremstyle{remark}

\author{Andrew Suk\thanks{Courant Institute, New York and EPFL, Lausanne. Email: {\tt suk@cims.nyu.edu}.  The author gratefully acknowledges the support from the Swiss National
Science Foundation, Grant No. 200021-125287/1.}
}

\title{A note on geometric 3-hypergraphs}

\begin{document}

\maketitle

\begin{abstract}

In this note, we prove several Tur\'an-type results on geometric hypergraphs.  The two main theorems are 1) Every $n$-vertex geometric 3-hypergraph in the plane with no three strongly crossing edges has at most $O(n^2)$ edges, 2) Every $n$-vertex geometric 3-hypergraph in 3-space with no two disjoint edges has at most $O(n^2)$ edges.  These results support two conjectures that were raised by Dey and Pach, and by Akiyama and Alon.

\end{abstract}

\section{Introduction}

A \emph{geometric $r$-hypergraph} $H$ \emph{in $d$-space} is a pair $(V,E)$, where $V$ is a set of points in general position in Euclidean $d$-space, and $E$ is a set of closed $(r-1)$-dimensional simplices (edges) induced by some $r$-tuple of $V$.  The sets $V$ and $E$ are called the \emph{vertex set} and \emph{edge set} of $H$, respectively.  Two edges in $H$ are \emph{crossing} if they are vertex disjoint and have a point in common.  Notice that if $k$ edges are pairwise crossing, it does not imply that they all have a point in common.  Hence we say that $H$ contains \emph{$k$ strongly crossing edges} if $H$ contains $k$ vertex disjoint edges that all share a point in common.  See Figure 1.

 \begin{figure}[h]
  \centering
  \label{fig:example}
    \subfigure[Three strongly crossing edges.]{\label{fig:tcase1}\includegraphics[width=.24\textwidth]{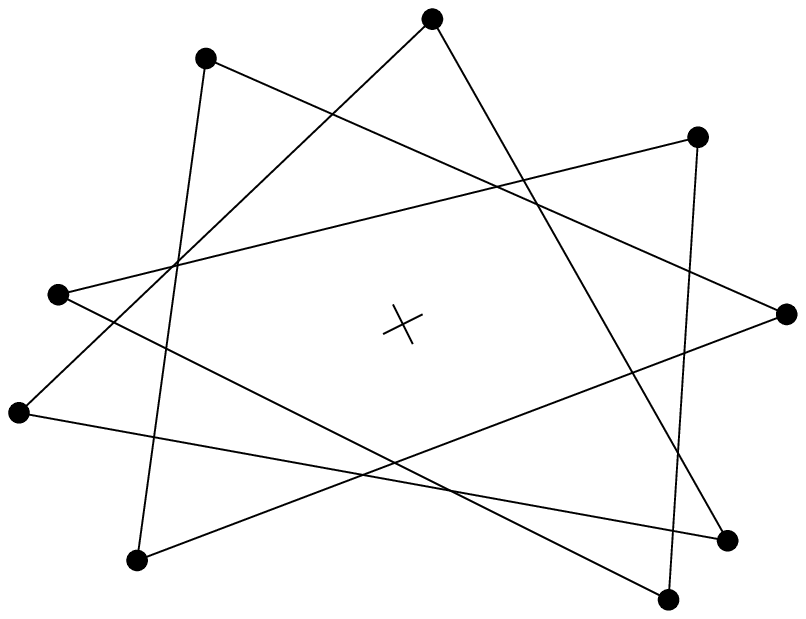}} \hspace{1cm}
        \subfigure[Three pairwise crossing edges with an empty intersection.]{\label{fig:tcase1}\includegraphics[width=.24\textwidth]{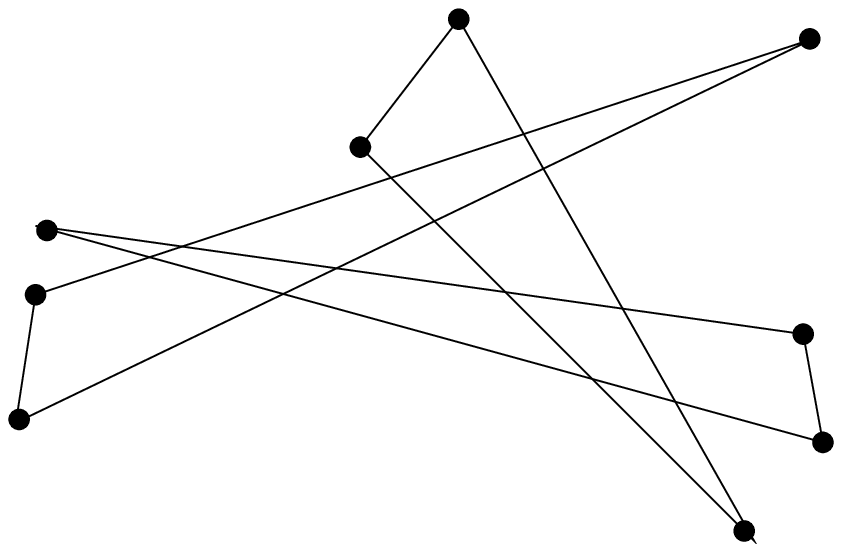}}\hspace{1cm}
    \subfigure[Three edges not strongly crossing since two share a vertex.]{\label{fig:tcase2}\includegraphics[width=.24\textwidth]{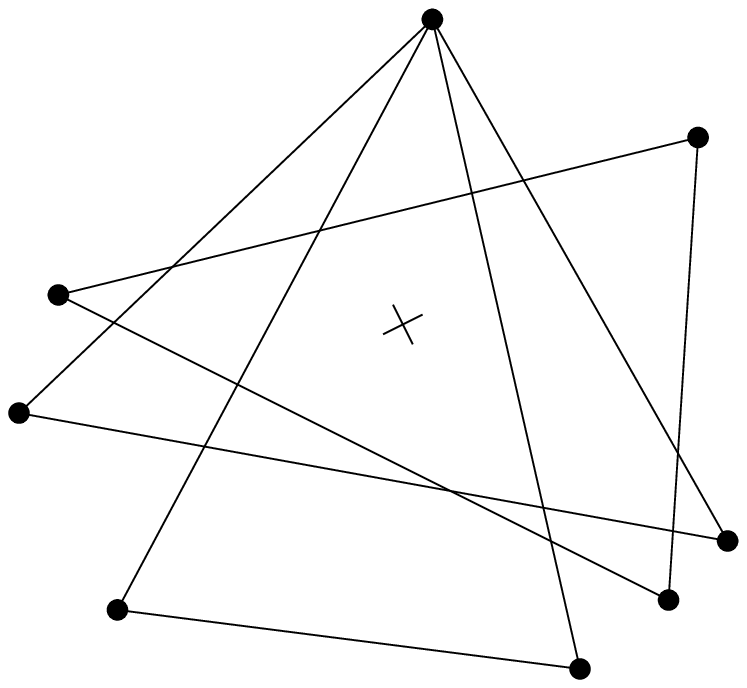}}
        \caption{Three edges of a geometric 3-hypergraph in the plane.}
\end{figure}

A direct application of the colored Tverberg theorem (see \cite{alon},\cite{z}) gives

\begin{theorem}
\label{deypach}
Let $ex_{d}(SC_k^{d+1},n)$ denote the maximum number of edges an $n$-vertex geometric $(d+1)$-hypergraph in $d$-space has with no $k$ strongly crossing edges.  Then

$$ex_{d}(SC_k^{d+1},n) = O\left(n^{d + 1 - \frac{1}{(2k - 1)^{d}}}\right).$$

\end{theorem}

\noindent Dey and Pach \cite{dey} showed that $ex_{d}(SC_2^{d+1},n) = \Theta(n^d)$, and conjectured $ex_{d}(SC_k^{d+1},n) = \Theta(n^{d})$ for every fixed $d$ and $k$.  The lower bound can easily be seen by taking all edges with a vertex in common.  The main motivation for their conjecture is for deriving upper bounds on the maximum number of \emph{$k$-sets} of an $n$-point set in $\mathbb{R}^d$.  See \cite{mat} for more details.  In this note, we settle the Dey-Pach conjecture for geometric 3-hypergraphs in the plane with no three strongly crossing edges, and improve the upper bound of $ex_2(SC_k^3,n)$.

\begin{theorem}
\label{first}
$ex_2(SC_3^3,n) = \Theta(n^2)$.

\end{theorem}

\begin{theorem}
\label{strong}
For fixed $k \geq 4$, $ex_{2}(SC_k^{3},n) \leq       O(n^{3 - \frac{1}{k }})$.

\end{theorem}

As a related result, Akiyama and Alon \cite{aki} used the Borsuk-Ullam Theorem \cite{bor} to show the following.

\begin{theorem}

Let $ex_d(D_k^d,n)$ denote the maximum edges that an $n$-vertex geometric $d$-hypergraph in $d$-space has with no $k$ pairwise disjoint edges.  Then

$$ex_d(D_k^d,n) \leq n^{d-(1/k)^{d-1}}.$$

\end{theorem}

\noindent They conjecture that for every fixed $d$ and $k$, $ex_d(D_k^d,n) = \Theta(n^{d-1})$.  Again the lower bound can easily be seen by taking all edges with a vertex in common.  Pach and T\"or\H ocsik \cite{pach} showed that $ex_2(D_k^2,n) = O(k^4n)$, which was later improved to $O(k^2n)$ by T\'oth \cite{geza}.  Here we settle the Akiyama-Alon conjecture for geometric 3-hypergraphs in 3-space with no two disjoint edges.

\begin{theorem}
\label{32}
 $ex_3(D_2^3,n)= \Theta(n^2).$

                          \end{theorem}

\noindent For clarity of the proofs, we do not make any attempts to optimize the constants.

\section{Strongly crossing edges in the plane}

In this section we will prove Theorems \ref{first} and \ref{strong}.  Recall that a \emph{geometric graph} is a graph drawn in the plane with vertices represented by points and edges by straight line segments connecting the corresponding pairs.  Recently Ackerman \cite{four} showed the following.

\begin{lemma}
\label{ack}

Let $G = (V,E)$ be an $n$-vertex geometric graph in the plane with no four pairwise crossing edges.  Then $|E(G)| \leq O(n)$.

\end{lemma}

$\hfill\square$

\begin{figure}
\begin{center}\label{level}
\includegraphics[width=320pt]{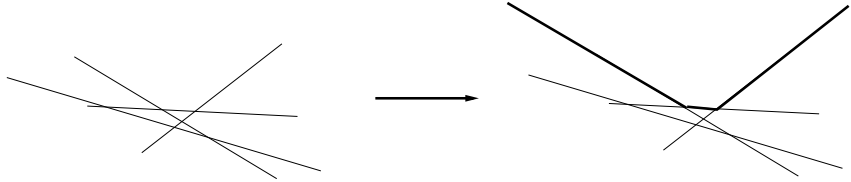}
  \caption{The top level of four pairwise crossing segments is drawn thick.}
 \end{center}
\end{figure}

\noindent We note that Lemma \ref{ack} holds for topological graphs.  Before we give the proofs, we will introduce some terminology.  Consider a family $\mathcal{S} = \{s_1,...,s_k\}$ of pairwise crossing segments in the plane, and let $\mathcal{L} = \{l_1,...,l_k\}$ be a family of lines such that $l_i$ is the line supported by segment $s_i$.  Recall that the \emph{level} of a point $x \in \cup \mathcal{L}$ is defined as the number of lines of $\mathcal{L}$ lying strictly below $x$.   We define the \emph{top level of} $\mathcal{L}$ as the closure of the set of points in $\cup \mathcal{L}$ with level $k-1$.  We define the \emph{top level of} $\mathcal{S}$ to be the top level of $\mathcal{L}$.  See Figure \ref{level}.  Notice that $L$ is a (not strictly) convex function.

\begin{figure}[h]
\begin{center}\label{base2}
\includegraphics[width=150pt]{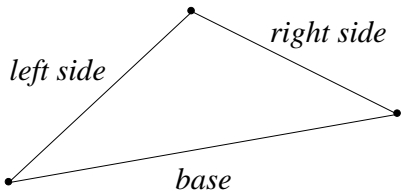}
  \caption{The base, left side, and right side.}
 \end{center}
\end{figure}

For each edge $t$ in a geometric 3-hypergraph in the plane, we define its \emph{base} as the side with the longest $x$-projection.  We define the other two sides of $t$ as its \emph{left} and \emph{right} side.  See Figure \ref{base2}.  Notice that every edge in a geometric 3-hypergraph is incident to a vertex that lies strictly above or below its base.  We are now ready to prove Theorem \ref{first}.

\medskip

\noindent \textbf{Proof of Theorem \ref{first}.}  Let $H = (V,E)$ be an $n$-vertex geometric $3$-hypergraph in the plane with no three strongly crossing edges.  We can assume that $|E(H)| \geq 20n^2$ (since otherwise we would be done) and at most $|E(H)|/2$ edges in $H$ are incident to a vertex that lies strictly below its base.  We will discard all such edges, leaving us with at least $|E(H)|/2$ edges left.  Let $E_{uv}$ be the set of edges in $H$ with base $uv$.  We discard all sets $E_{uv}$ for which $|E_{uv}| \leq |E(H)|/(2n^2)$.  Since we have thrown away at most $|E(H)|/4$ edges in this process, we have at least $|E(H)|/4$ edges left.  Therefore $|E_{uv}| = 0$ or $|E_{uv}| \geq |E(H)|/(2n^2) \geq 10$.

Now let $G_v = (V,E)$ denote the geometric graph with $V(G_v) = V(H)$ and $xy \in E(G_v)$ if $conv(x\cup y\cup v) \in E(H)$ with base $xy$.

 \begin{observation}
\label{obs}
 $G_v$ does not contain four pairwise crossing edges (bases).

 \end{observation}

\noindent \textbf{Proof.}  For sake of contradiction, suppose $G_v$ contains four pairwise crossing edges $b_1,b_2,b_3,b_4\in E(G_v)$.   Then $v$ lies above $b_i$ for all $i$.  Let $L$ denote the top level of the arrangement $\mathcal{S} = \{b_1,b_2,b_3,b_4\}$.  Now the proof falls into three cases.
\medskip

\noindent \emph{Case 1.}  Suppose $L$ intersects exactly two members of $\mathcal{S}$, say bases $b_1$ and $b_2$ (in order from left to right along $L$).  Let $p$ be the intersection point of $b_1$ and $b_2$.  Then the vertical line through $p$ must intersect $b_3$ below $p$.  Moreover, since segments $b_1$ and $b_3$ cross, $v$ and the right-endpoint of $b_3$ must lie on the same half-plane generated by the line supported by $b_1$.  Likewise, $v$ and the left-endpoint of $b_3$ must lie on the same half-plane generated by the line supported by $b_2$.  Therefore  $p \in conv(v\cup b_3)$.  See Figure \ref{fig:tcase1}.  Since $|E_{b_1}|, |E_{b_2}| \geq 10$, there exists vertices $x,y\in V(H)$ such that $conv(v\cup b_3), conv(x\cup b_1), conv(y\cup b_2)$ are three (vertex disjoint) strongly crossing edges in $H$ and we have a contradiction.

\medskip

\noindent \emph{Case 2.}  Suppose $L$ intersects exactly three members of $\mathcal{S}$, say bases $b_1,b_2,b_3$ (in order from left to right along $L$).  Now $b_4$ must intersect $b_2$ either to the left or right of $b_2\cap L$.  Without loss of generality, we can assume that $b_4$ intersects $b_2$ to the right of $b_2\cap L$.  Let $p$ be the intersection point of segments $b_2$ and $b_3$.  By the same argument as above,  $p \in conv(v\cup b_4)$.  See Figure \ref{fig:tcase2}.  Since $|E_{b_1}|, |E_{b_2}| \geq 10$, there exists vertices $x,y\in V(H)$ such that $conv(v\cup b_4), conv(x\cup b_1), conv(y\cup b_2)$ are three strongly crossing edges in $H$ and we have a contradiction.

\medskip

\begin{figure}[h]
  \centering
  \label{fig:threepics}
    \subfigure[Case 1.]{\label{fig:tcase1}\includegraphics[width=.3\textwidth]{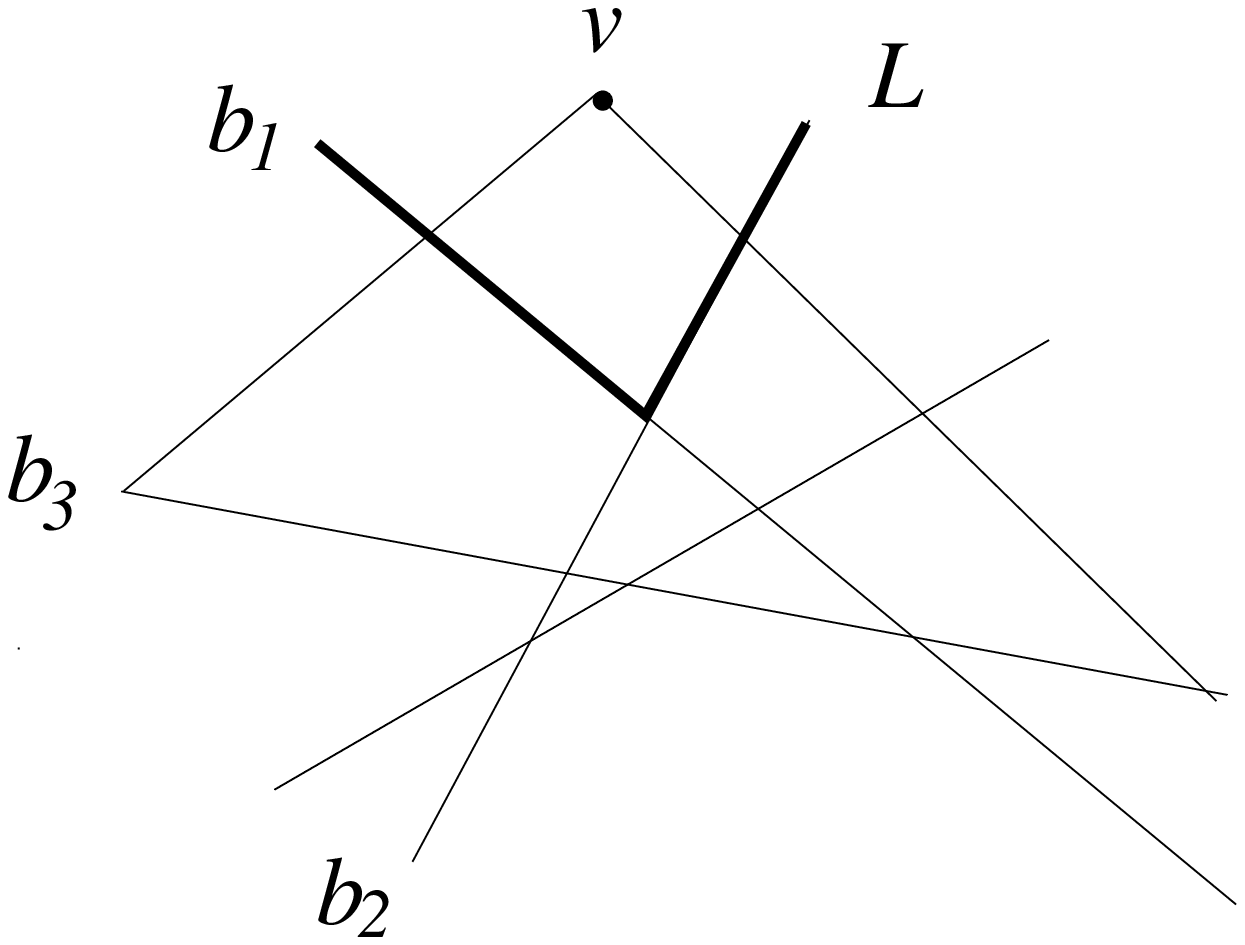}} \hspace{.5cm}
    \subfigure[Case 2.]{\label{fig:tcase2}\includegraphics[width=.3\textwidth]{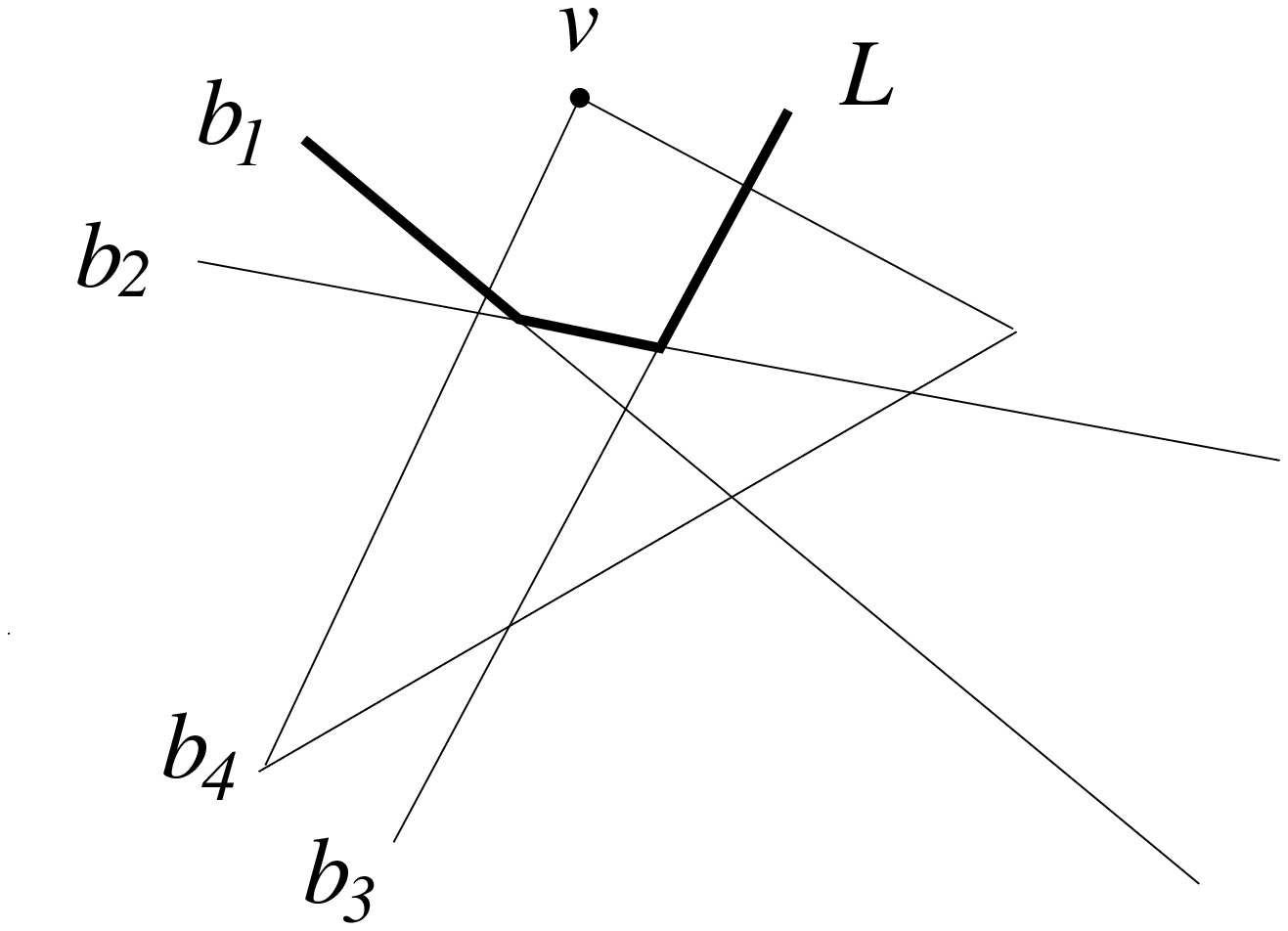}}   \hspace{.5cm}
        \subfigure[Case 3.]{\label{fig:tcase3}\includegraphics[width=.3\textwidth]{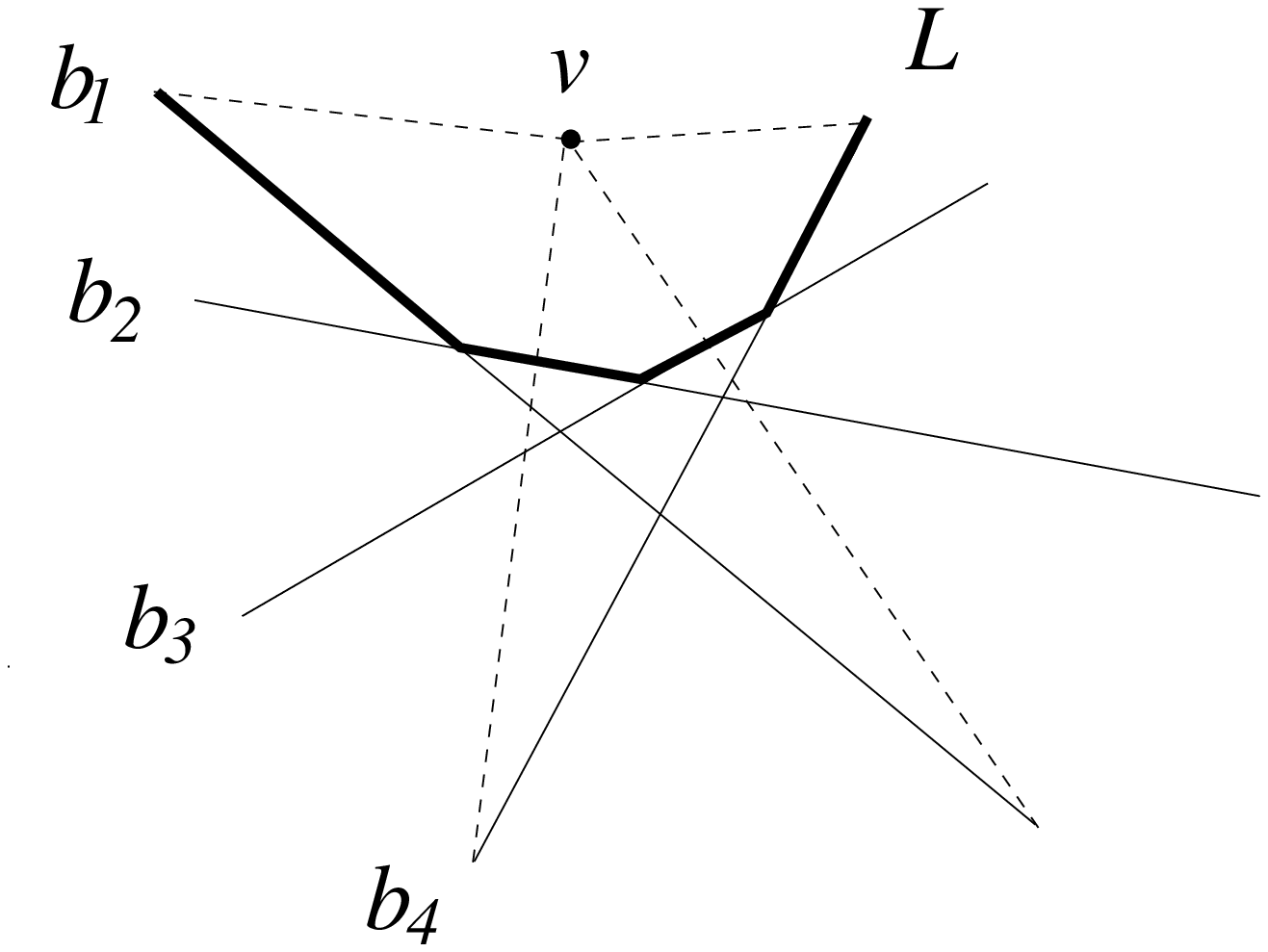}}
        \caption{Three cases.}
\end{figure}

\noindent \emph{Case 3.}   Suppose $L$ intersects $b_1,b_2,b_3,b_4$ in order from left to right along $L$.  Let $p$ be the intersection point of segments $b_2$ and $b_3$, and let $l$ be the vertical line through $v$.  Since the right endpoint of $b_4$ lies to the right of $l$, and the left endpoint of $b_1$ lies to the left of $l$, we have $p \in conv(v\cup b_1)\cup conv(v\cup b_4)$.  Therefore, either $conv(v\cup b_1)$ or $conv(v\cup b_4)$ (say $conv(v\cup b_1)$) contains $p$.  See Figure \ref{fig:tcase3}.  Since $|E_{b_2}|, |E_{b_3}| \geq 10$, there exists vertices $x,y\in V(H)$ such that $conv(v\cup b_1), conv(x\cup b_2), conv(y\cup b_3)$ are three strongly crossing edges in $H$ and we have a contradiction.

\medskip

$\hfill\square$

\noindent Therefore by Lemma \ref{ack}, $|E(G_v)| \leq O(n)$ for every vertex $v \in V(H)$.  Hence

$$\frac{|E(H)|}{4} \leq \sum\limits_{v\in V(H)} |E(G_v)| = O(n^2),$$

\noindent which implies $|E(H)| = O(n^2)$.

$\hfill\square$

\medskip

Before we prove Theorem \ref{strong}, we will need the following lemma due to Valtr \cite{valtr}.

\begin{lemma}
\label{valtr}

Let $G = (V,E)$ be an $n$-vertex geometric graph in the plane such that all of the edges in $G$ intersect the $y$-axis.  If $G$ does not contain $k$ pairwise crossing edges, then $|E(G)| \leq c_k n$ where $c_k$ depends only on $k$.

\end{lemma}

$\hfill\square$

\medskip

\noindent \textbf{Proof of Theorem \ref{strong}}. Let $H$ be an $n$-vertex geometric 3-hypergraph in the plane with no $k$ strongly crossing edges for $k\geq 4$.  Just as before, we can assume at most $|E(H)|/2$ of the edges in $H$ are incident to a vertex that lies strictly below its base.  We discard all such edges, leaving us with at least $|E(H)|/2$ edges left in $H$.  Now we make the following observation.

\begin{observation}
\label{key}
Suppose $b_1,...,b_k$ are $k$ pairwise crossing bases and $v_1,...,v_k \in V(H)$ such that $conv(v_i\cup b_j) \in E(H)$ with base $b_j$ for all $i,j$.  Then $H$ contains $k$ strongly crossing edges.

\end{observation}

\begin{figure}[h]
  \centering
\includegraphics[width=0.4\textwidth]{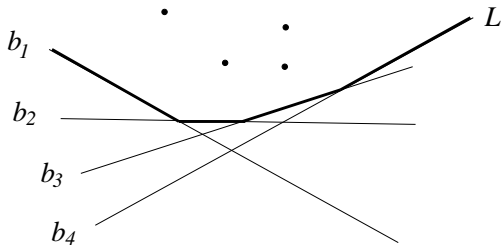}
    \caption{Arrangement of $b_1,b_2,b_3,b_4$.}
        \label{fig:level}
\end{figure}

\noindent \textbf{Proof.}  Let $L$ denote the top level of the segment arrangement $\mathcal{S} = \{b_1,...,b_k\}$ and assume that $b_1,...,b_k$ are ordered by increasing slopes.  See Figure \ref{fig:level}.

Now we define edges $t_1,t_2,...,t_k\in E(H)$ as follows.  Among the $k$ edges $conv(b_1\cup v_1), conv(b_1\cup v_2),..., conv(b_1\cup v_k) \in E(H)$, (with slight abuse of notation) let $t_1 = conv(b_1\cup v_1)$ be the edge whose right side has the rightmost intersection with $L$.  Then among the $k-1$ edges $conv(b_2\cup v_2), conv(b_2\cup v_3),..., conv(b_2\cup v_k)$, (again with slight abuse of notation) let $t_2 = conv(b_2\cup v_2)$ be the edge whose right side has the rightmost intersection with $L$.  We continue this procedure until we have $k$ edges $t_1,t_2,...,t_k$.  Clearly these $k$ edges are vertex disjoint.

\begin{figure}[h]
  \centering
    \includegraphics[width=0.85\textwidth]{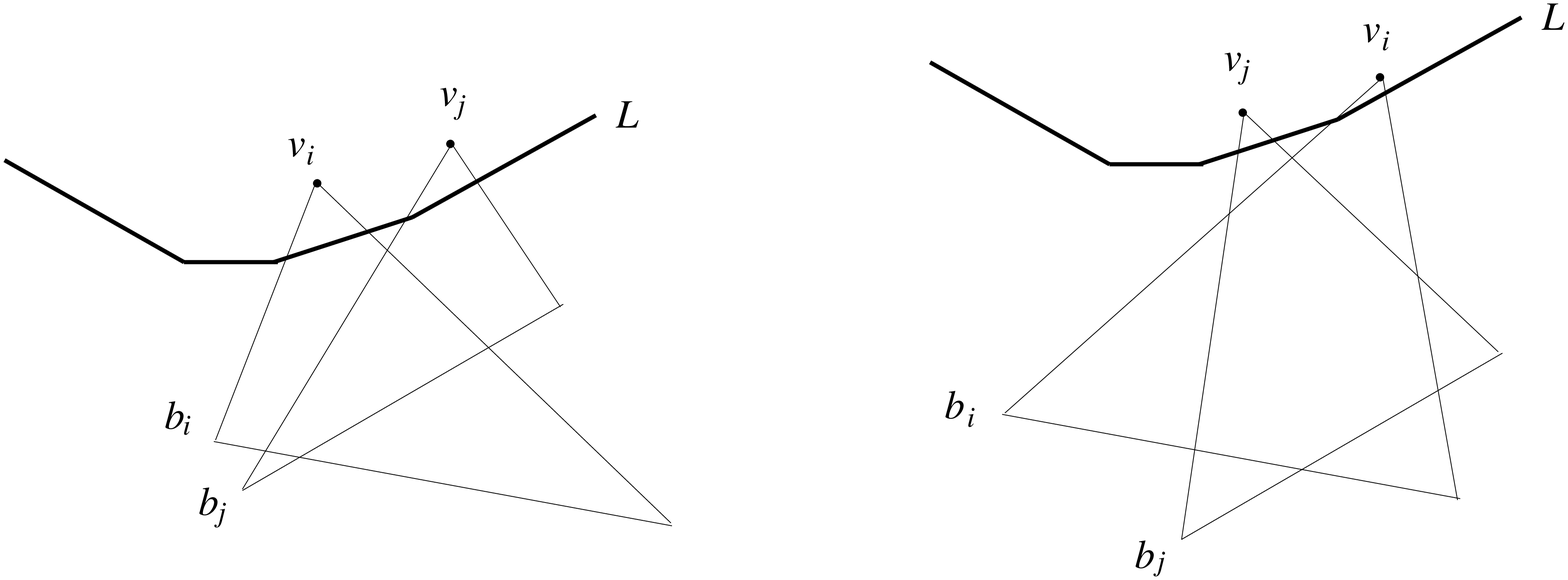}
    \caption{Assume $(t_i\cap L)\cap(t_j\cap L) = \emptyset$.}
    \label{fig:miss}
\end{figure}

Now notice that $(t_i\cap L)\cap(t_j\cap L) \neq \emptyset$ for all pairs $i,j$.  Indeed for sake of contradiction, suppose there exists two edges $t_i$ and $t_j$ for $i < j$ such that either $t_i\cap L$ lies completely to the left of $t_j\cap L$ or vice versa.  See Figure \ref{fig:miss}.

\medskip

\noindent \emph{Case 1.}   Suppose $t_i\cap L$ lies completely to the left of $t_j\cap L$.  Then the vertical line through $v_j$ intersects the right side of $t_i$ below $v_j$. Therefore the right side of $conv(b_i\cup v_j)$ intersects $L$ more to the right than the right side of $t_i = conv(b_i\cup v_i)$ does.  This contradicts the definition of $t_i$ and $t_j$.

\medskip

\noindent \emph{Case 2.}  Suppose $t_i\cap L$ lies completely to the right of $t_j\cap L$.  Then there exists a base $b_s$ that has a point $p$ on $L$ between $t_i\cap L$ and $t_j\cap L$.  Base $b_s$ must

  \begin{enumerate}

  \item lie below $v_i$ and $v_j$,

  \item cross $b_i$ and $b_j$, and

  \item contain point $p$.

  \end{enumerate}

  However this impossible by the following argument.  Let $l$ be the vertical line through $p$.  Clearly $l$ intersects $b_i$ and $b_j$.  Since $b_s$ lies below $v_i$ and $v_j$, $b_s$ must intersect $b_j$ to the left of $l$, and intersect $b_i$ to the right of $l$.  Since $b_s$ intersects $b_j$ to the left of $l$, the slope of $b_s$ must be greater than the slope of $b_j$.  However since the slope of $b_i$ is less than the slope of $b_j$, this implies that $b_s$ cannot intersect $b_i$ to the right of $l$.  Hence we have a contradiction.

\begin{figure}[h]
  \centering
    \includegraphics[width=0.3\textwidth]{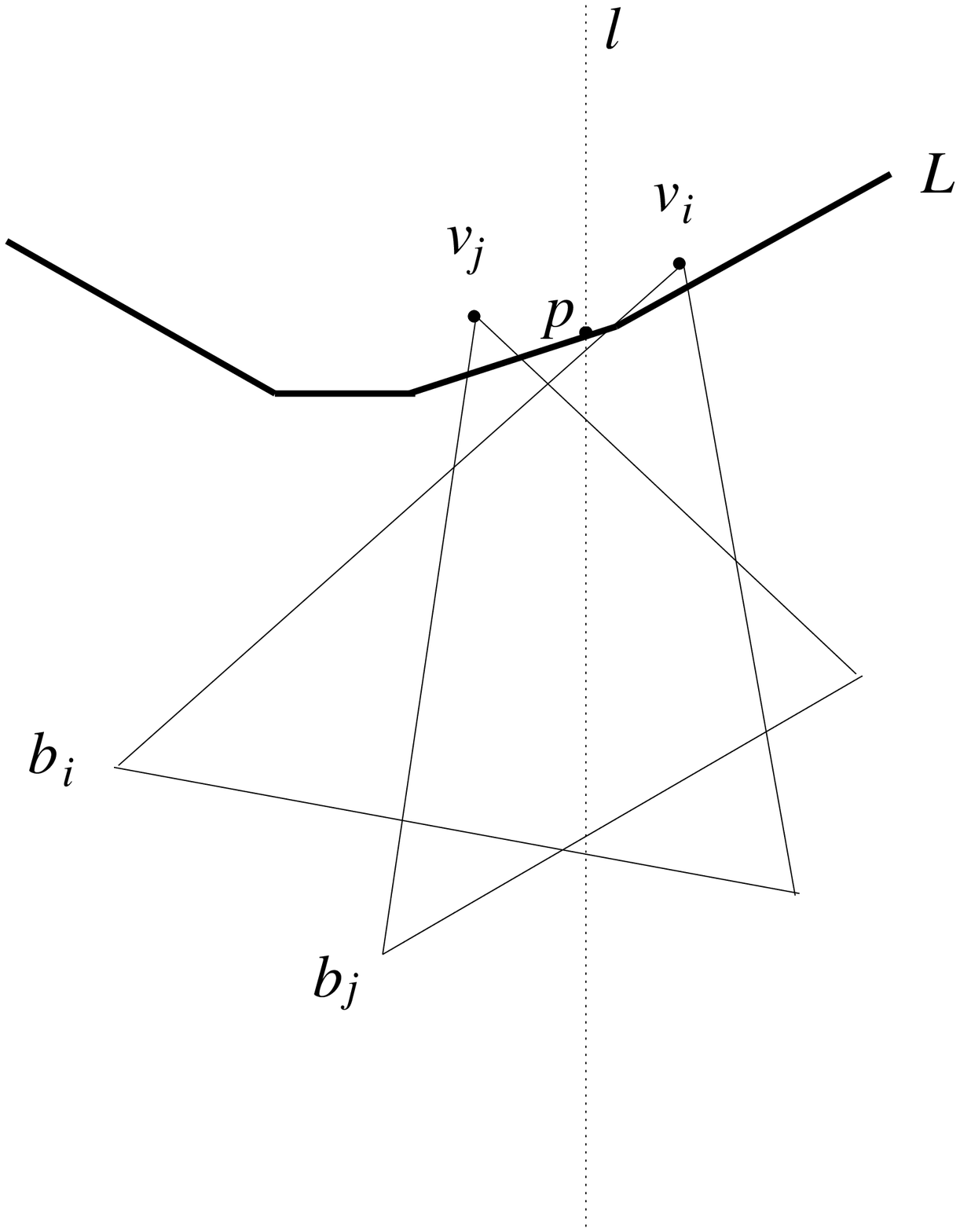}
    \caption{Case 2.}
    \label{fig:counterb}
\end{figure}

\medskip

\noindent Since $(t_i\cap L) \cap (t_j\cap L) \neq \emptyset$ for every $i,j \in \{1,2,..,k\}$, by Helly's Theorem \cite{helly} $t_1,...,t_k$ has a nonempty intersection on $L$.

$\hfill\square$

Notice that no $k$ points in $V(H)$ have $c_kn$ bases in common. Indeed, otherwise the vertical line through any of these $k$ points would intersect all $c_kn$ bases, and by Lemma \ref{valtr} there would be $k$ pairwise crossing bases.  By Observation \ref{key}, we would have $k$ strongly crossing edges.

Now let $G = (A\cup B,E)$ be a bipartite graph where $A = V(H)$ and $B = V^2(H)$, such that $(v,xy) \in E(G)$ if $conv(x\cup y\cup v) \in E(H)$ with base $xy$.  Since $G$ does not contain $K_{k,c_kn}$ as a subgraph, we can use the following well known result of K\H ov\'ari, S\'os, Tur\'an \cite{sos}.

\begin{theorem}
\label{kovari}

If $G = (A\cup B,E)$ is a bipartite graph with $|A| = n$ and $|B| = m$ containing no subgraph $K_{r,s}$ with the $r$ vertices in $A$ and the $s$ vertices in $B$, then

$$|E(G)| \leq (s-1)^{1/r}nm^{1- 1/r} + (r-1)m.$$

\end{theorem}

\noindent By plugging in the values $m = n^2, r = k, s = c_kn$ into Theorem \ref{kovari}, we obtain

$$\frac{|E(H)|}{2} \leq |E(G)| \leq O\left(n^{3-\frac{1}{k}}\right).$$

\noindent Hence

$$|E(H)|\leq O\left(n^{3-\frac{1}{k}}\right).$$

$\hfill\square$

\subsection{Convex geometric 3-hypergraphs}

In the case when the vertices are in convex position in the plane, extremal problems on geometric 3-hypergraphs become easier due to the linear ordering of its vertices.  The proof of Observation \ref{key} can be copied almost verbatim to conclude the following.

\begin{observation}
\label{kcross}
Let $H  = (V,E)$ be a geometric 3-hypergraph in the plane with vertices in convex position.  Suppose $H$ contains $k$ edges of the form $t_i = conv(x_i\cup y_i \cup z_i)$, such that the vertices $(x_1,...,x_k,y_1,...,y_k,z_1,...,z_k)$ appear in clockwise order along the boundary of their convex hull.   Then $t_1,...,t_k$ are $k$ strongly crossing edges.

\end{observation}

$\hfill\square$

Marcus and Klazar \cite{klazar} extended the Marcus-Tardos theorem \cite{tardos} by showing that the number of 1-entries in a $r$-dimensional $(0,1)$-matrix with side length $n$ which avoids an $r$-dimensional permutation matrix is $O(n^{r-1})$.  As pointed out by Marcus and Klazar, it is not difficult to modify their proof to obtain an $O(n^{r-1})$ bound on the number of edges in an ordered $n$-vertex $r$-uniform hypergraph that does not contain a fixed ordered matching. Hence by Observation \ref{kcross}, we can conclude the following.

\begin{theorem}
\label{convex}
Let $H = (V,E)$ be a geometric 3-hypergraph in the plane with vertices in convex position.  If $H$ does not contain $k$ strongly crossing edges, then $|E(H)| \leq c_k n^2$ where $c_k$ is a constant that depends only on $k$.

\end{theorem}

$\hfill\square$

\section{Disjoint edges in 3-space}

In this section, we will prove Theorem \ref{32}.   Recall that two edges in a geometric graph are \emph{parallel} if they are the opposite edges of a convex quadrilateral.  Katchalski and Last  \cite{last} and Pinchasi \cite{pin} showed that all $n$-vertex geometric graphs with more than $2n-2$ edges contain two parallel edges.  By following Pinchasi's argument almost verbatim, one can prove the following.

\begin{lemma}
\label{par}
Let $G$ be a graph drawn on the unit sphere $S$ with vertices represented as points such that no three lie on a great circle, and edges $uv \in E(G)$ are drawn as arcs along the great circle containing points $u$ and $v$ of length less than $\pi$ (the shorter arc).  We say that edges $e_1,e_2\in E(G)$ are avoiding if the great circle supported by $e_1$ is disjoint to $e_2$, and the great circle supported by $e_2$ is disjoint from $e_1$.  If $|E(G)|  > 2n-2$, then $G$ contains two avoiding edges.

\end{lemma}
$\hfill\square$

\noindent \textbf{Proof of Theorem \ref{32}.}  Let $H = (V,E)$ be an $n$-vertex geometric 3-hypergraph in 3-space with no two disjoint edges.  Fix a pair of vertices $u,v \in V(H)$, and just consider the edges $E_{uv} =\{t\in E(H): \textnormal{$u,v$ are vertices of $t$}\}$.  We color $t \in E_{uv}$ \emph{red} if all of the members of $E_{uv}$ lie in one of the closed half-spaces generated by the plane supported by $t$.  Notice that there are at most two red edges in $E_{uv}$.  Repeat this procedure for each pair of vertices, which will leave us with at most $n^2$ red edges in the end.  Color the remaining edges blue, and let $d_b(v)$ denote the number of blue edges incident to $v$.  Then we have

$$\sum\limits_{v \in V(H)} d_b(v) \geq 3E(H) - 3n^2.$$

\noindent Therefore, there exists a vertex $v$ incident to at least $(3|E(H)| - 3n^2)/n$ blue edges.  Now consider a small 2-dimensional sphere $S^2$ centered at $v$.  Then the intersection of $S^2$ and the blue edges incident to $v$ forms a graph $G$ with at most $n$ vertices and at least $(3E(H) - 3n^2)/n$ edges.

\begin{figure}[h]
  \centering
\includegraphics[width=0.25\textwidth]{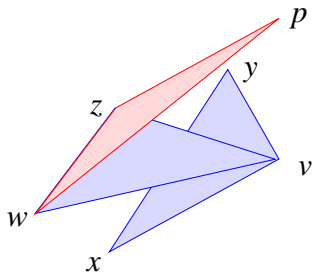}
    \caption{Disjoint edges $conv(w\cup z\cup p)$ and $conv(x\cup y\cup v)$.}
        \label{fig:disjoint}
\end{figure}

If $(3|E(H)| - 3n^2)/n > 2n-2$, then by Lemma \ref{par} we know that $G$ contains two avoiding edges $xy$ and $wz$.  Let $h$ be the plane supported by the blue edge $conv(w\cup z\cup v) \in E(H)$.  Then the blue edge $conv(x\cup y\cup v)$ must lie in one of the closed half-spaces generated by the plane $h$.  Since $conv(w\cup z\cup v)$ is blue, there must be a red edge $conv(w\cup z\cup p)$ such that $h$ separates it from $conv(x\cup y\cup v)$. Hence $conv(x\cup y\cup v)$ and $conv(w\cup z\cup p)$ are disjoint and we have a contradiction.  See Figure \ref{fig:disjoint}.  Therefore $(3|E(H)| - 3n^2)/n \leq 2n-2$, which implies $|E(H)| \leq O(n^2)$.

$\hfill\square$

\section{Remarks}

By applying the Abstract Crossing Lemma (see \cite{wagner}) to Theorem \ref{first}, every $n$-vertex geometric 3-hypergraph $H$ in the plane has either $O(n^2)$ edges or $\Omega(|E(H)|^7/n^{12})$ triples that have a point in common.  In the latter case, by the fractional Helly theorem \cite{kat} this implies one can always find a point inside at least $\Omega(|E(H)|^5/n^{12})$ edges of $H$.  However, this is not as strong as the

$$\Omega\left(\frac{|E(H)|^3}{n^{6}\log^2n}\right)$$

\noindent bound obtained by Nivasch and Sharir \cite{niv}.

\end{document}